%% Plain TeX, with the macro package xypic

\newcount\secno
\newcount\prmno
\newif\ifnotfound
\newif\iffound

\def\namedef#1{\expandafter\def\csname #1\endcsname}
\def\nameuse#1{\csname #1\endcsname}

\long\def\ifundefined#1#2#3{\expandafter\ifx\csname
  #1\endcsname\relax#2\else#3\fi}
\def\hwrite#1#2{{\let\the=0\edef\next{\write#1{#2}}\next}}

% Working with lists
\toksdef\ta=0 \toksdef\tb=2
\long\def\leftappenditem#1\to#2{\ta={\\{#1}}\tb=\expandafter{#2}%
                                \edef#2{\the\ta\the\tb}}
\long\def\rightappenditem#1\to#2{\ta={\\{#1}}\tb=\expandafter{#2}%
                                \edef#2{\the\tb\the\ta}}

\def\lop#1\to#2{\expandafter\lopoff#1\lopoff#1#2}
\long\def\lopoff\\#1#2\lopoff#3#4{\def#4{#1}\def#3{#2}}

\def\ismember#1\of#2{\foundfalse{\let\given=#1%
    \def\\##1{\def\next{##1}%
    \ifx\next\given{\global\foundtrue}\fi}#2}}

% Les commandes
\def\section#1{\medbreak
               \global\def\currenvir{section}
               \global\advance\secno by1\global\prmno=0
               {\bf \number\secno. {#1}}
               \smallskip}

\def\subsection{\global\def\currenvir{subsection}
                \global\advance\prmno by1
                \ind{ (\number\secno.\number\prmno) }}
\def\subsec{\global\def\currenvir{subsection}
                \global\advance\prmno by1
                { (\number\secno.\number\prmno)\ }}

\def\proclaim#1{\global\advance\prmno by 1
                {\bf #1 \the\secno.\the\prmno$.-$ }}

\long\def\th#1 \enonce#2\endth{%
   \medbreak\proclaim{#1}{\it #2}\global\def\currenvir{th}\smallskip}

\def\rem#1{\global\advance\prmno by 1
{\it #1} \the\secno.\the\prmno$.-$ }

% CROSS-REFERENCES

\def\isinlabellist#1\of#2{\notfoundtrue%
   {\def\given{#1}%
    \def\\##1{\def\next{##1}%
    \lop\next\to\za\lop\next\to\zb%
    \ifx\za\given{\zb\global\notfoundfalse}\fi}#2}%
    \ifnotfound{\immediate\write16%
                 {Warning - [Page \the\pageno] {#1} No reference found}}%
                \fi}%
\def\ref#1{\ifx\labellist\empty{\immediate\write16
                 {Warning - No references found at all.}}
               \else{\isinlabellist{#1}\of\labellist}\fi}

\def\newlabel#1#2{\rightappenditem{\\{#1}\\{#2}}\to\labellist}
\def\labellist{}

\def\label#1{\relax}

 \def\openall{\openout\lbl=\jobname.lbl}

\newread\testfile
\def\lookatfile#1{\openin\testfile=\jobname.#1
    \ifeof\testfile{\immediate\openout\nameuse{#1}\jobname.#1
                    \write\nameuse{#1}{}
                    \immediate\closeout\nameuse{#1}}\fi%
    \immediate\closein\testfile}%

\def\begin{\newlabel{theta}{1.2}
\newlabel{Pic}{1.3}
\newlabel{lin}{1.4}
\newlabel{Trd}{2.1}
\newlabel{Tr0}{2.2}
\newlabel{normal}{2.3}
\newlabel{D_r}{2.4}
\newlabel{Verlinde}{3.1}
\newlabel{sum}{3.2}
\newlabel{formule}{3.3}
\newlabel{total}{3.4}
\newlabel{rem}{4.3}}

 %% Fin de la numerotation automatique

\magnification 1250
\mathcode`A="7041 \mathcode`B="7042 \mathcode`C="7043
\mathcode`D="7044 \mathcode`E="7045 \mathcode`F="7046
\mathcode`G="7047 \mathcode`H="7048 \mathcode`I="7049
\mathcode`J="704A \mathcode`K="704B \mathcode`L="704C
\mathcode`M="704D \mathcode`N="704E \mathcode`O="704F
\mathcode`P="7050 \mathcode`Q="7051 \mathcode`R="7052
\mathcode`S="7053 \mathcode`T="7054 \mathcode`U="7055
\mathcode`V="7056 \mathcode`W="7057 \mathcode`X="7058
\mathcode`Y="7059 \mathcode`Z="705A
\def\spacedmath#1{\def\packedmath##1${\bgroup\mathsurround =0pt##1\egroup$}
\mathsurround#1
\everymath={\packedmath}\everydisplay={\mathsurround=0pt}}
\def\nospacedmath{\mathsurround=0pt
\everymath={}\everydisplay={} } \spacedmath{2pt}
\font\eightrm=cmr8         \font\eighti=cmmi8
\font\eightsy=cmsy8        \font\eightbf=cmbx8
        \font\eightit=cmti8
        \font\sixrm=cmr6
\font\sixi=cmmi6           \font\sixsy=cmsy6
\font\sixbf=cmbx6
\catcode`\@=11
\def\eightpoint{%
  \textfont0=\eightrm \scriptfont0=\sixrm \scriptscriptfont0=\fiverm
  \def\rm{\fam\z@\eightrm}%
  \textfont1=\eighti  \scriptfont1=\sixi  \scriptscriptfont1=\fivei
    \textfont2=\eightsy \scriptfont2=\sixsy \scriptscriptfont2=\fivesy
   \textfont\itfam=\eightit
  \def\it{\fam\itfam\eightit}%
   \textfont\bffam=\eightbf \scriptfont\bffam=\sixbf
  \scriptscriptfont\bffam=\fivebf
  \def\bf{\fam\bffam\eightbf}%
   \abovedisplayskip=9pt plus 3pt minus 9pt
  \belowdisplayskip=\abovedisplayskip
  \abovedisplayshortskip=0pt plus 3pt
  \belowdisplayshortskip=3pt plus 3pt
  \smallskipamount=2pt plus 1pt minus 1pt
  \medskipamount=4pt plus 2pt minus 1pt
  \bigskipamount=9pt plus 3pt minus 3pt
  \normalbaselineskip=9pt
  \setbox\strutbox=\hbox{\vrule height7pt depth2pt width0pt}%
  \let\bigf@nt=\eightrm     \let\smallf@nt=\sixrm
  \normalbaselines\rm}

\def\note#1#2{\footnote{\parindent
0.4cm$^#1$}{\vtop{\eightpoint\baselineskip12pt\hsize15.5truecm\noindent #2}}
\parindent 0cm}
\def\sdir_#1^#2{\mathrel{\mathop{\kern0pt\oplus}\limits_{#1}^{#2}}}
\def\pc#1{\tenrm#1\sevenrm}
\def\up#1{\raise 1ex\hbox{\smallf@nt#1}}
\def\tx{\kern-1.5pt -}
\def\cqfd{\kern 2truemm\unskip\penalty 500\vrule height 4pt depth 0pt width
4pt\medbreak} \def\carre{\vrule height 4pt depth 0pt width 4pt}
\def\no{n\up{o}\kern 2pt}
\def\ind{\par\hskip 1truecm\relax}

\def\rond{\kern 1pt{\scriptstyle\circ}\kern 1pt}
\def\iso{\mathrel{\mathop{\kern 0pt\longrightarrow }\limits^{\sim}}}

\def\det{\mathop{\rm det}\nolimits}
\def\Pic{\mathop{\rm Pic}\nolimits}

\def\dim{\mathop{\rm dim}\nolimits}

\def\Tr{\mathop{\rm Tr}\nolimits}

\def\ch{\mathop{\rm ch}\nolimits}

\font\gragrec=cmmib10
\def\Mu{\hbox{\gragrec \char22}}
\font\san=cmssdc10
\def\ext {\hbox{\san \char3}}
\def\sy{\hbox{\san \char83}}

\input amssym.def
\input amssym
\vsize = 25truecm
\hsize = 16truecm
%\hoffset = -.15truecm
\voffset = -.5truecm
\parindent=0cm
\baselineskip15pt
\overfullrule=0pt
\begin
\centerline{\bf The Verlinde formula for ${\bf PGL}_p$}
\smallskip
\smallskip \centerline{Arnaud {\pc BEAUVILLE\note{1}{Partially supported by the
European HCM project ``Algebraic Geometry in Europe" (AGE).}}} \vskip0.9cm
\hfill\vbox{\hsize4cm\eightpoint\baselineskip12pt\line{\hfill\it To the memory
of\hfill} \line{\hfill Claude ITZYKSON\hfill}}
\vskip0.7cm

{\bf Introduction}
\smallskip

\ind The Verlinde formula expresses the number of linearly independent
conformal
blocks in any rational conformal field theory. I am concerned here with a quite
particular case, the Wess-Zumino-Witten model associated to a complex
semi-simple
group\note{2}{This group is the complexification of the compact semi-simple
group
considered by physicists.} $G$. In this case the space of conformal blocks can
be
interpreted as the space of holomorphic sections of a line bundle on a
particular
projective variety, the moduli space $M_G$ of holomorphic $G$\tx bundles on the
given Riemann surface.  The fact that the dimension of this space of sections
can be explicitly computed is of great interest for mathematicians, and a
number of
rigorous proofs of that formula (usually called by mathematicians, somewhat
incorrectly, the ``Verlinde formula") have been  recently given (see e.g.\ [F],
[B-L], [L-S]).  \ind These proofs deal only with simply-connected groups.
In this paper we treat the case of  the projective  group
${\bf PGL}_r$ when $r$ is prime.   \ind Our approach is to relate to the case
of ${\bf SL}_r$, using standard algebro-geometric methods. The
components $M_{{\bf PGL}_r}^d$ $(0\le  d<r)$ of the moduli space $M_{{\bf
PGL}_r}$  can be
identified with the quotients $M_r^d/J_r$ , where  $M_r^d$ is the moduli space
of
vector bundles on $X$ of rank $r$ and fixed determinant of degree $d$, and
$J_r$ the
finite group of holomorphic line bundles $\alpha $ on $X$ such that $\alpha
^{\otimes
r}$ is trivial. The space we are looking for is the space of $J_r$\tx invariant
global sections of a line bundle ${\cal L}$ on $M_r^d$; its dimension can be
expressed in terms of the  character of the representation of $J_r$ on
 $H^0(M_r^d,{\cal L})$. This is given by the Lefschetz trace formula, with a
subtlety for $d=0$, since $M_r^0$ is not smooth.
The key point  (already used in [N-R]) which makes the computation quite easy
is that
the fixed point set of any non-zero element of $J_r$ is an abelian variety  --
this is
where the assumption on the group is essential. Extending the method to other
cases
would require a Chern classes computation on the moduli space  $M_H$ for some
semi-simple subgroups $H$ of $G$; this may be feasible, but goes far beyond the
scope
of the present  paper. Note that the case of $M_{{\bf PGL}_2}^1$
has been previously worked out in [P] (with an unfortunate misprint in the
formula).

 \ind In the last section we check that our formulas agree with the
predictions of Conformal Field Theory, as they appear for instance in [F-S-S].
Note
that our results are slightly more precise (in this  particular case): we get a
formula  for $\dim H^0(M_{{\bf PGL}_r}^d,{\cal L})$ for every $d$, while   CFT
only
predicts the sum of these dimensions (see Remark \ref{rem}).

\vskip1cm \section{The moduli space $M_{{\bf PGL}_r}$} \subsection Throughout
the
paper we denote by $X$  a compact (connected) Riemann surface, of genus $g\ge
2$; we
fix a point $p$ of $X$.  Principal ${\bf PGL}_r$\tx bundles on $X$ correspond
in a
one-to-one way to projective bundles of rank $r-1$ on $X$, i.e. bundles of the
form
${\bf P}(E)$, where $E$ is a rank $r$ vector bundle on $X$; we say that  ${\bf
P}(E)$ is semi-stable if the vector bundle $E$ is semi-stable. The semi-stable
projective  bundles of rank $r-1$ on $X$ are parameterized by a projective
variety,   the moduli space $M_{{\bf PGL}_r}$. \ind   Two vector bundles $E$,
$F$
give rise to isomorphic projective bundles if and only if $F$ is isomorphic to
$E\otimes \alpha $ for some line bundle $\alpha $ on $X$. Thus a projective
bundle
can always be written as ${\bf P}(E)$ with  $\det E={\cal O}_X(dp)$,  $0\le
d<r$;  the
vector bundle $E$ is then determined up to tensor product by a line bundle
$\alpha $
with $\alpha ^r={\cal O}_X$. In particular, the moduli space $M_{{\bf
PGL}_r}$ has $r$ connected components $M_{{\bf PGL}_r}^d$ $(0\le  d<r)$. Let us
denote
by $M_r^d$  the moduli space of semi-stable vector bundles on $X$ of rank $r$
and
determinant ${\cal O}_X(dp)$, and by $J_r$  the kernel of the multiplication by
$r$
in the Jacobian $JX$ of $X$; it is a finite group, canonically isomorphic to
$H^1(X,{\bf Z}/(r))$. The group $J_r$ acts on $M_r^d$, by the rule $(\alpha
,E)\mapsto E\otimes\alpha $; it follows from the above remarks that the
component
$M_{{\bf PGL}_r}^d$ is isomorphic to the quotient
$M_r^d/J_r$. \subsection\label{theta} We will need a precise description of the
line
bundles on $M_{{\bf PGL}_r}$. Let me first recall how one describes line
bundles on
$M_r^d$ [D-N]: a simple way  is to mimic the classical
definition of the theta divisor on the Jacobian of $X$ (i.e.\ in the rank $1$
case).
Put $\delta=(r,d)$; let $A$ be a  vector bundle on $X$ of rank $r/\delta$ and
degree $(r(g-1)-d)/\delta$. These conditions imply $\chi (E\otimes A)=0$ for
all $E$
in $M_r^d$; if $A$ is general enough, it follows that the condition
$H^0(X,E\otimes
A)\not=0$ defines a (Cartier) divisor $\Theta _A$ in $M_r^d$. The corresponding
line
bundle ${\cal L}_d:={\cal O}(\Theta _A)$ does not depend on the choice of $A$,
and
 generates the Picard group $\Pic(M_r^d)$.
\subsection\label{Pic}
The quotient map $q:M_r^d\rightarrow M_{{\bf PGL}_r}^d $ induces a
homomorphism\break $q^*:\Pic(M_{{\bf PGL}_r}^d)\rightarrow \Pic(M_r^d)$,
which is easily seen to be injective. Its image is determined in [B-L-S]: it is
generated by ${\cal L}_d^\delta$ if $r$ is odd, by ${\cal L}_d^{2\delta}$ if
$r$ is
even.

\subsection\label{lin} Let ${\cal L}'$ be a line bundle on $M_{{\bf
PGL}_r}^d$. The  line
bundle ${\cal L}:=q^*{\cal L}'$ on $M_r^d$ admits a natural  action of
$J_r$, compatible with the action of $J_r$ on $M_r^d$ (this is often
called a $J_r$\tx linearization of ${\cal L}$). This action is characterized by
the
property that every element $\alpha $ of $J_r$ acts trivially on the fibre of
${\cal
L}$ at a point of $M_r^d$ fixed by $\alpha $. In the sequel we will always
consider
 line bundles  on $M_r^d$  of the form $q^*{\cal L}'$, and endow them
with the above $J_r$\tx linearization.

\ind This linearization defines a representation of $J_r$ on the space of
global sections;
essentially by definition, the global sections of ${\cal L}'$
correspond to the $J_r$\tx invariant sections of ${\cal L}$. Therefore our task
will be to compute the dimension of the space of invariant sections; as
indicated
in the introduction, we will do that by computing, for any $\alpha \in J_r$ of
order $r$, the trace of $\alpha $ acting on  $H^0(M_r^d,{\cal L})$.

\vskip1cm
 \section {The action of
$J_r$ on $H^0(M_r^d,{\cal L}_d^k)$}

\ind We start with the case when $r$ and $d$ are coprime, which is easier to
deal with because the moduli space is smooth.
\th Proposition
\enonce Assume $r$ and $d$ are coprime. Let $k$ be an integer; if $r$ is even
we
assume that $k$ is even. Let $\alpha $ be an element of order $r$ in $JX$. Then
the
trace of $\alpha $ acting on $H^0(M_r^d,{\cal L}_d^{k})\ $ is
$(k+1)^{(r-1)(g-1)}$.
\endth\label{Trd}
{\it Proof}:   The
Lefschetz trace formula reads [A-S]$$\Tr(\alpha\, |\,H^0(M_r^d,{\cal
L}_d^k))=\int_{P}{\rm Todd}(T_P)\ \lambda (N_{P/M_r^d},\alpha )^{-1}\
\widetilde{\ch}({\cal L}^k_{d\,|P},\alpha )\ .$$
 Here $P$ is the fixed subvariety
of $\alpha $;  whenever $F$ is a vector bundle on $P$ and $\varphi $ a
diagonalizable endomorphism of $F$, so that $F$ is the direct sum of its
eigen-sub-bundles $F_\lambda $ for $\lambda \in {\bf C}$, we put
  $$\widetilde{\ch}(F,\varphi  )=\sum \lambda \,{\ch}(F_\lambda )\quad ;\quad
\lambda (F,\varphi )=\prod_\lambda \sum_{p\ge   0} (-\lambda )^p
\ch(\ext^pF_\lambda ^*)\
.$$
\ind We have a number of informations on the right hand side thanks to
[N-R]:

(\ref{Trd} {\it a})
\ Let
$\pi :\widetilde{X}\rightarrow X$ be the \'etale $r$\tx sheeted covering
associated to
$\alpha $; put $\xi =\alpha ^{r(r-1)/2}\in JX$.  The map $L\mapsto \pi _*(L)$
identifies any component  of the fibre of the
norm map  ${\rm Nm}:J^d\widetilde{X}\rightarrow J^dX$ over $\xi(dp)$ with $P$.
In particular, $P$ is isomorphic to an abelian variety, hence the term ${\rm
Todd}(T_P)$ is trivial.

(\ref{Trd} {\it b}) \  Let $\theta \in H^2(P,{\bf Z})$ be the restriction
to $P$ of the class of the principal polarization of $J^d\widetilde{X}$. The
term
$\lambda (N_{P/M_r^d},\alpha )$ is equal to $r^{r(g-1)}e^{-r\theta }$.

(\ref{Trd} {\it c}) \ The dimension of $P$ is $N=(r-1)(g-1)$, and one has
$\int_P
{\theta^N \over N!}=r^{g-1}$.

\ind  With our convention the action of $\alpha $ on ${\cal L}^k_{d\,|P}$ is
trivial.
The  class $c_1({\cal L}_{d\,|P})$ is equal to $r\theta $: the pull back to $P$
of the
theta divisor $\Theta _A$ (\ref{theta}) is the divisor of line bundles $L$ in
$P$
with $H^0(L\otimes \pi ^*A)\not=0$; to compute its cohomology class we may
replace
$\pi ^*A$ by any vector bundle with the same rank  and degree,  in particular
by a
direct sum of $r$ line bundles of degree $r(g-1)-d$, which gives the required
formula.
\ind Putting things together,  we find
$$\Tr(\alpha\, |\,H^0(M_r^d,{\cal L}_d^k))=\int_{P} r^{-r(g-1)}e^{r\theta }
e^{kr\theta } = (k+1)^{(r-1)(g-1)}\ .\quad  \carre$$

\bigskip
\ind We now consider the degree $0$ case:

\th Proposition
\enonce Let  $k$ be a multiple of $r$, and of $2r$ if $r$ is even; let $\alpha
$ be
an element of order $r$ in $JX$. Then the trace of $\alpha $ acting on
$H^0(M_r^0,{\cal L}_0^{k})\ $ is $ ({k\over r}+1)^{(r-1)(g-1)}$.
\endth\label{Tr0}
{\it Proof}: We cannot
apply directly the Lefschetz trace formula since it is manageable
only for smooth projective varieties; instead we use another well-known tool,
the
Hecke correspondence (this idea appears for instance in [B-S]). For simplicity
we
write $M_d$ instead of $M_r^d$.  There
exists a Poincar\'e bundle ${\cal E}$ on $X\times M_1$, i.e.\ a vector bundle
whose
restriction to $X\times\{E\}$, for each point $E$ of $M_1$, is isomorphic to
$E$.
Such a  bundle is determined up to tensor product by a line bundle coming from
$M_1$; we will see later how to normalize it. We denote by ${\cal E}_p$ the
restriction of ${\cal E}$ to $\{p\}\times M_1$, and by
  ${\cal P}$ the projective
bundle  ${\bf P}({\cal E}_p^*)$ on $M_1$. A point of ${\cal P}$ is a pair
$(E,\varphi )$ where $E$ is a vector bundle in $M_1$ and $\varphi :E\rightarrow
{\bf C}_p$ a non-zero homomorphism, defined up to a scalar; the kernel of
$\varphi $
is then a vector bundle $F\in M_1$, and we can view equivalently a point of
${\cal P}$
as a pair of vector bundles $(F,E)$ with $F\in M_0$, $E\in M_1$ and $F\i E$.
The
projections $p_d$ on $M_d$ $(d=0,1)$ give rise to the ``Hecke diagram"
\input xypic $$\diagram
&{\cal P}\dlto_{p_1} \drto^{p_0} &&\\    M_1 &&M_0&\kern-40pt.\\  \enddiagram$$

\smallskip
\th Lemma
\enonce  The Poincar\'e bundle ${\cal E}$ can be normalized (in a unique way)
so
that $\det {\cal E}_p={\cal L}_1\ ;$ then ${\cal O}_{\cal P}(1) \cong
p_0^*{\cal
L}_0$. \endth \label{normal}
{\it Proof}:  Let $E\in M_1$. The fibre $p_1^{-1}(E)$ is the projective space
of
non-zero linear forms $\ell :E_p\rightarrow {\bf C}$, up to a scalar.
The restriction of $p_0^*{\cal L}_0$ to this projective space  is ${\cal O}(1)$
(choose a  line bundle $L$ of degree $g-1$ on $X$; if $E$ is general enough,
$H^0(X,E\otimes L)$ is spanned by a section $s$ with $s(p)\not=0$, and the
condition
that the bundle $F$ corresponding   to $\ell $ belongs to $\Theta _L$ is the
vanishing of $\ell (s(p))$). Therefore $p_0^*{\cal L}_0$ is of the form
${\cal O}_{\cal P}(1)\otimes p_1^*{\cal N}$ for some line bundle ${\cal N}$ on
$M_1$.
Replacing ${\cal E}$ by ${\cal E}\otimes{\cal N}$ we ensure ${\cal O}_{\cal
P}(1)
\cong p_0^*{\cal L}_0$.
\ind An easy computation gives
 $K_{\cal  P}=p_1^*{\cal  L}_1^{-1}\otimes p_0^*{\cal  L}_0^{-r}$ ([B-L-S],
Lemma 10.3).  On the other hand, since ${\cal P}={\bf P}({\cal E}_p^*)$, one
has
$K_{\cal  P}= p_1 ^*(K_{M_1}\otimes
\det {\cal  E}_p)\otimes {\cal  O}_{\cal  P}(-r)$;  using $K_{
M_1}={\cal  L}^{-2}_1$ [D-N], we get $\det{\cal  E}_p={\cal  L}_1$. \cqfd

\medskip
\ind We normalize ${\cal E}$ as in the lemma; this gives for each $k\ge   0$ a
canonical
isomorphism $p_{1*}p_0^*{\cal L}_0^k\cong \sy^k{\cal E}_p$. Let $\alpha $ be an
element of  order $r$ of  $JX$. It acts on the various
moduli spaces in sight; with a slight abuse of language, I will still denote by
$\alpha $ the corresponding automorphism.  There exists an isomorphism  $\alpha
^*{\cal E}\iso {\cal E}\otimes\alpha $, unique up to a scalar ([N-R], lemma
4.7);
 the induced isomorphism $u:\alpha ^*{\cal
E}_p\iso {\cal E}_p$ induces the action of $\alpha $ on ${\cal
P}$. Imposing $u^r={\rm Id}$ determines $u$ up to a $r$\tx th root of unity,
hence
determines completely $\sy^ku$ when $k$ is a multiple of $r$.   Since the Hecke
diagram is equivariant with respect to $\alpha $, it gives rise to a diagram of
isomorphisms
 $$\diagram  & H^0({\cal P},p_0^*{\cal L}^k_0) & \\
  H^0(M_1,\sy^k{\cal E}_p) \urto^{p_1^*}  & & H^0(
  M_0,{\cal L}_0^k)\ulto_{p_0^*}
  \enddiagram$$
which is compatible with the action of $\alpha $; in particular, the trace we
are
looking for is equal to the trace of $\alpha $ on  $H^0(M_1,\sy^k{\cal E}_p)$.
\ind We are now in the situation of Prop.\ \ref{Trd},
and the Lefschetz trace formula gives:
$$\Tr(\alpha\, |\,H^0(M_1,\sy^k{\cal E}_p))=\int_{P}{\rm
Todd}(T_P)\ \lambda (N_{P/M_1},\alpha )^{-1}\
\widetilde{\ch}(\sy^k{\cal E}_{p\,|P},\alpha )\ .$$
\ind  The only term we  need to compute is $\widetilde{\ch}(\sy^k{\cal
E}_{p\,|P},\alpha )$. Let ${\cal N}$ be the restriction to $\widetilde{X}\times
P$
of a Poincar\'e line bundle on $\widetilde{X}\times J^1\widetilde{X}$; let us
still
denote by  $\pi :
\widetilde{X}\times P\rightarrow  X\times P$ the map $\pi \times {\rm Id}_P$.
The
vector bundles $\pi _*({\cal N})$ and ${\cal E}_{|X\times P}$ have the same
restriction to $ X\times \{\gamma \}$ for all $\gamma \in P$, hence after
tensoring
${\cal N}$ by
 a line bundle on $P$ we may assume they are isomorphic ([R], lemma 2.5).
Restricting
to $\{p\} \times P$ we get $\displaystyle {\cal E}_{p\,|P}=\sdir_{\pi
(q)=p}^{}{\cal
N}_{q}$, with ${\cal N}_q={\cal N}_{|\{q\}\times P}$. \smallskip

\ind We claim that the ${\cal N}_q$'s are the eigen-sub-bundles of ${\cal
E}_{p\,|P}$
relative to $\alpha $. By (\ref{Trd} {\it a}), a pair $(E,F)\in{\cal P}$ is
fixed by
$\alpha $ if and only if $E=\pi_*L$, $F=\pi_*L'$, with ${\rm Nm}(L)=\xi(p) $,
${\rm
Nm}(L')=\xi $; because of the inclusion $F\i E$ we may take $L'$ of the form
$L(-q)$, for some point $q\in \pi^{-1}(p)$. In other words, the fixed locus of
$\alpha $ acting on ${\cal P}$ is the disjoint union of the sections $(\sigma
_q)_{q\in \pi^{-1}(p)}$ of the fibration $p_1^{-1}(P)\rightarrow P$
characterized by
$\sigma _q(\pi_*L)=(\pi_*L,\pi_*(L(-q)))$. Viewing ${\cal
P}$ as ${\bf P}({\cal
E}^*_{p\,|P})$, the section $\sigma _q$ corresponds to
the exact sequence $$0\rightarrow \pi _*({\cal N}(-q))_{\,|\{p\}\times
P}\longrightarrow
 \pi _*({\cal N})_{\,| \{p\}\times P}\cong {\cal E}_{\,|\{p\}\times
P}\longrightarrow {\cal N}_q\rightarrow 0\ .$$
 Therefore on each fibre ${\bf P}(E_p)$, for
$E\in P$, the automorphism $\alpha $ has exactly $r$ fixed points,
corresponding to
the $r$ sub-spaces ${\cal
N}_{(q,E)}$ for $q\in \pi ^{-1} (p)$; this proves our claim.
\ind  The line bundles ${\cal N}_q$ for $q\in \widetilde{X}$ are  algebraically
equivalent, and therefore have the same Chern class. We thus have $c_1(  {\cal
E}_{p\,|P})=r\,c_1({\cal N}_q)$. On the other hand we know that
 $\det
{\cal E}_{p}={\cal L}_1$ (lemma \ref{normal}), and that $c_1({\cal
L}_{1\,|P})=r\theta $ (proof of Prop.\ \ref{Trd}). By comparison we get
$c_1({\cal N}_q)=\theta $. Putting things together we obtain
$$\widetilde{\ch}(\sy^k{\cal E}_{p\,|P},\alpha )=\int_P \Tr \sy^kD_r\
e^{k\theta
}r^{-r(g-1)}e^{r\theta }$$  where $D_r$ is the  diagonal $r$\tx by-$\!r$ matrix
with entries  the $r$ distinct $r$\tx th roots of unity.

\th Lemma
\enonce  The trace of
$\sy^kD_r$ is $1$ if $r$ divides $k$ and $0$ otherwise.
\endth\label{D_r}
\ind Consider the formal series $\displaystyle s(T):=\sum_{i\ge
0}T^i\Tr\sy^iu$ and
$\displaystyle \lambda (T):=\sum_{i\ge   0} T^i\Tr\ext^iu$. The formula
$s(T)\lambda
(-T)=1$ is well-known (see e.g. [Bo], \S 9, formula (11)). But  $$\lambda
(-T)=\sum_{i=0}^r(-T)^i\Tr\ext^iu=\prod_{\zeta^r=1}(1-\zeta T)=1-T^r\ ,$$
hence the lemma.  Using (\ref{Trd} {\it c}) the Proposition follows. \cqfd

\vskip1cm
\section{Formulas}
\ind In this section I will apply the above results to compute the dimension of
the space of sections of the line bundle ${\cal L}_d^k$ on the moduli space
$M_{{\bf
PGL}_r}^d$.  Let me first recall the corresponding Verlinde   formula for the
moduli
spaces
$M_r^d$. Let $\delta=(r,d)$; we write ${\cal L}_d={\cal D}^{r/\delta}$, with
the
convention that we only consider powers of ${\cal D}$ which are multiple of
$r/\delta$ (the line bundle ${\cal D}$ actually makes sense on the {\it moduli
stack}
${\cal M}_r^d$, and generates its Picard group). We denote by $\Mu_r$ the
center of
${\bf SL}_r$, i.e.\ the group of scalar matrices $\zeta {\it I}_r$ with $\zeta
^r=1$. \th Proposition \enonce Let $T_k$ be the set of diagonal matrices
 $t={\rm diag}(t_1,\ldots,t_r)$ in ${\bf SL}_r({\bf C})
$ with  $t_i\not=t_j$ for $i\not=j$, and
$t^{k+r}\in \Mu_r$; for $t\in T_k$, let
$\displaystyle \delta(t)=\prod_{i<j}(t_i-t_j)$. Then
$$\dim H^0(M_r^d,{\cal D}^k)=r^{g-1}(k+r)^{(r-1)(g-1)}\  \sum_{
t\in T_k/{\goth S}_r} { ((-1)^{r-1} t^{k+r})^{-d} \over
|\delta(t)|^{2g-2}}\quad
\cdot$$ \endth\label{Verlinde}
{\it Proof}: According to  [B-L],
Thm. 9.1, the space $H^0(M_r^d,{\cal D}^k)$ for $0<d<r$ is canonically
isomorphic to
the space of conformal blocks in genus $g$ with  the
 representation $V_{k\varpi_{r-d}}$ of ${\bf SL}_r$ with highest weight
$k\varpi_{r-d}$ inserted at one point. The Verlinde formula gives therefore
(see [B],
Cor. 9.8\note{1}{There is a misprint  in the first equality of that corollary,
where one should read $T_\ell ^{\rm reg}/W$ instead of $T_\ell ^{\rm reg}$;
the
second equality (and the proof!) are correct.}\kern-4pt): $$\dim
H^0(M_r^d,{\cal
D}^k)=r^{g-1}(k+r)^{(r-1)(g-1)}\sum_{t\in T_k/{\goth S}_r }
{\Tr^{}_{V_{k\varpi_{r-d}}}(t)\over |\delta(t)|^{2g-2}}\ ;$$ this is still
valid
for $d=0$ with the convention $\varpi_r=0$. \ind  The character of the
representation $V_{k\varpi_{r-d}}$ is given by the Schur formula (see e.g.
[F-H],
Thm. 6.3): $$\Tr^{}_{V_{k\varpi_{r-d}}}(t)={1\over \delta (t)}\ \left|
\matrix{t_1^{k+r-1} & t_2^{k+r-1} & \ldots & t_r^{k+r-1}\cr t_1^{k+r-2}
&t_2^{k+r-2}&\ldots & t_r^{k+r-2}\cr \vdots & \vdots & & \vdots \cr t_1^{k+d}
&t_2^{k+d}&\ldots & t_r^{k+d}\cr t_1^{d-1} &t_2^{d-1}&\ldots & t_r^{d-1}\cr
\vdots &
\vdots & & \vdots \cr 1 &1&\ldots & 1\cr
}
\right|\ .$$
\ind Writing $t^{k+r}=\zeta{\it I}_r \in\Mu_r$, the big determinant
reduces to \break $\zeta ^{r-d}(-1)^{d(r-d)}\det(t_j^{d-i})$, and finally,
since
$\prod t_i=1$, to  $((-1)^{r-1}\zeta) ^{-d}\delta (t)$, which gives the
required formula. \cqfd \smallskip
\th Corollary
\enonce Let $T'_k$ be the set of matrices $t={\rm diag}(t_1,\ldots,t_r)$ in
${\bf SL}_r({\bf C})$ with   $t_i\not=t_j$ if $i\not=j$, and
$t^{k+r}=(-1)^{r-1}{\it I}_r$. Then $$\sum_{d=0}^{r-1}\dim H^0(M_r^d,{\cal
D}^k)=r^{g}(k+r)^{(r-1)(g-1)}\ \sum_{t\in T'_k/{\goth S}_r} { 1 \over
|\delta(t)|^{2g-2}}\ \cdot\quad \carre$$
 \endth \label{sum}\smallskip

\ind We now consider the  moduli space $M_{{\bf PGL}_r}$. We know that
the line bundle ${\cal D}^k$ on $M_r^d$ descends to $M_{{\bf
PGL}_r}^d=M_r^d/J_r$
exactly when $k$ is a multiple of $r$ if $r$ is odd, or of $2r$ if $r$ is even
(\ref{Pic}). When this is the case we obtain a line bundle on $M_{{\bf
PGL}_r}^d$,
that we will still denote by ${\cal D}^k$; its global sections correspond to
the
$J_r$\tx invariant sections of $H^0(M_r^d,{\cal D}^{k})$.
 \ind We will assume that $r$ is {\it
prime}, so that
 every non-zero element $\alpha $ of $J_r$ has order $r$. Then Prop.\ \ref{Trd}
and
\ref{Tr0} lead immediately to a formula for the dimension  of
the $J_r$\tx invariant subspace of $H^0(M_r^d,{\cal D}^{k})$ as the average of
the
numbers $\Tr(\alpha )$ for $\alpha $ in $J_r$. Using Prop.\ \ref{Verlinde} we
conclude:  \th Proposition
\enonce  Assume that $r$ is prime. Let  $k$ be  a
multiple of $r$; if $r=2$ assume
 $4\mid k$. Then
$$\nospacedmath\displaylines{\dim H^0(M_{{\bf PGL}_r}^d,{\cal D}^{k})
 =  r^{-2g}\,\dim H^0(M_r^d,{\cal D}^k)\ +\ (1-r^{-2g})({k\over
r}+1)^{(r-1)(g-1)}\cr
\hfill = r^{-2g}\,({k\over r}+1)^{(r-1)(g-1)}\ \Bigl(r^{r(g-1)}
 \sum_{ t\in T_k/{\goth S}_r} { ((-1)^{r-1}t^{k+r})^{-d} \over
|\delta(t)|^{2g-2}} \, +\, r^{2g}-1 \Bigr)\ .}$$
\endth\label{formule}
\smallskip
\ind Summing over $d$  and plugging in Cor.\ \ref{sum} gives the following
rather
complicated formula:
\th Corollary
\enonce
$$\dim H^0(M_{{\bf PGL}_r},{\cal D}^{k})=r^{1-2g}\,({k\over r}+1)^{(r-1)(g-1)}\
\Bigl( r^{
r(g-1)}\sum_{t\in T'_k/{\goth S}_r} { 1 \over
|\delta(t)|^{2g-2}}\ +\ r^{2g}-1\Bigr)\ .$$
\endth\label{total}
\ind  As an example, if  $k$ is an integer
divisible by $4$, we get $$\dim H^0(M_{{\bf PGL}_2},{\cal
D}^{k})=2^{1-2g}\,({k\over
2}+1)^{g-1}\bigl( \sum_{l\ {\rm odd}\atop 0<l<k+2} { 1 \over (\sin{l\pi \over
k+2})^{2g-2}}\,+\,2^{2g}-1\bigr)\ .\leqno(3.5)$$
\vskip1cm
\section{Relations with Conformal Field Theory}
 \subsection According to Conformal Field Theory, the space $H^0(M_{{\bf
PGL}_r},{\cal
D}^{k})$ should be canonically isomorphic to the space of conformal blocks for
a
certain Conformal Field Theory, the WZW model associated to the projective
group. This implies in particular that its dimension should be equal to
$\sum_j|S_{0j}|^{2-2g}$, where $(S_{ij})$ is a unitary symmetric matrix.
For instance in the case of the WZW model associated to ${\bf SL}_2$,
one has
$$S_{0j}={\sin {(j+1)\pi \over k+2}\over
\sqrt{{k\over 2}+1}}\qquad\hbox{, with}\quad 0\le  j\le  k\ , $$
where the index $j$ can be thought as running through the set of irreducible
representations $\sy^1,\ldots,\sy^k$ of ${\bf SL}_2$ (or equivalently ${\bf
SU}_2$),
with $\sy^j:=\sy^j({\bf C}^2)$.
\ind  We deduce from (3.5) an analogous expression for
${\bf PGL}_2$: we restrict ourselves to  even indices and write
$$S'_{0j}=2\,S_{0j}\qquad {\rm
for}\quad  j\hbox { even }<k/2\quad  {\rm ;}\qquad  S'_{0,{k\over
2}^{(1)}}=S'_{0,{k\over 2}^{(2)}} =S_{0{k\over 2}}\quad
\cdot$$ In other words,  we consider
only those representations of ${\bf SL}_2$ which factor through ${\bf PGL}_2$
and we
identify the representation $\sy^{2j}$ with $\sy^{k-2j}$,  doubling the
coefficient
$S_{0j}$ when these two representations are distinct, and counting twice the
representation
which is fixed by the involution (this process is  well-known, see
e.g.\ [M-S]).
 \subsection The case of ${\bf SL}_r$ is completely analogous; we only
need a few more terminology from representation theory (we follow the notation
of
[B]).  The primary fields are indexed by
 the set $P_k$ of dominant weights $\lambda $ with $\lambda
(H_\theta )\le  k$,   where $H_\theta $ is the matrix $\ {\rm
diag}(1,0,\ldots,0,-1)$.
For $\lambda \in P_k$, we put $\displaystyle t_\lambda =\exp 2\pi
i\,{\lambda+\rho\over k+r}$ (we identify the Cartan algebra of diagonal
matrices
with its dual using the standard bilinear form); the map $\lambda \mapsto
t_\lambda $ induces a bijection of $P_k$ onto $T_k/{\goth S}_r$ ([B], lemma 9.3
{\it c})). In view of Prop.\  \ref{Verlinde},  the coefficient
$S_{0\lambda }$ for $\lambda \in P_k$ is given by
$$S_{0\lambda } = {\delta(t_\lambda )\over
\sqrt{r}(k+r)^{(r-1)/2}}\ \cdot$$
\ind Passing to ${\bf PGL}_r$, we first restrict the indices to the subset
$P'_k$ of
 elements $\lambda \in P_k$  such that $t_\lambda $ belongs to $T'_k$;
this means that $\lambda $ belongs to the root lattice, i.e.\ that the
representation $V_\lambda$ factors through ${\bf PGL}_r$.
 The center $\Mu_r$ acts on $T_k$ by multiplication; this action preserves
$T'_k$,
and commutes with the action of ${\goth S}_r$. The corresponding action on
$P_k$ is
deduced, via the bijection $\lambda
\mapsto {\lambda +\rho \over k+r}$, from the standard action of $\Mu_r$ on the
fundamental alcove $A$ with vertices
$\{0,\varpi_1,\ldots,\varpi_{r-1}\}$.\note{1}{The element $\exp \varpi_1$ of
the
center gives the rotation of $A$ which maps $0$ to $\varpi_1$, $\varpi_1$ to
$\varpi_2$, $\ldots,$ and $\varpi_{r-1}$ to $0$.}

\ind We identify two elements of $P'_k$ if they are in the same orbit with
respect
to  this action. The action has a unique fixed point, the weight ${k\over
r}\rho
$, which corresponds to the diagonal matrix $D_r$ (\ref{D_r}); we associate to
this
weight
$r$ indices $\nu ^{(1)},\ldots,\nu ^{(r)}$, and put $$S'_{0\lambda
}=r\,S_{0\lambda
}\quad {\rm for}\   \lambda \in P'_k/\Mu_r\ ,\ \lambda \not={k\over r}\rho\
{\rm ;} \qquad  S'_{0,\nu ^{(i)} }= S_{0,{k\over r}\rho}\quad {\rm for}\
i=1,\ldots,r\ .$$
One deduces easily from Cor.\
\ref{total} the formula $\dim H^0(M_{{\bf PGL}_r},{\cal
D}^{k})=\sum |S'_{0\lambda }|^{2-2g}$, where $\lambda $ runs
over $P'_k/\Mu_r\cup\{\nu ^{(1)},\ldots,\nu ^{(r)}\}$.
\bigskip

\rem {Remark}\label{rem}
It is not clear to me what is the physical meaning of the space $H^0(M_{{\bf
PGL}_r}^d,{\cal D}^{k})$, in particular if its dimension can be predicted
in terms of the $S$\tx matrix. It is interesting to observe that the number
$N(g)$
given by Prop.\ \ref{formule}, which is equal to $\dim H^0(M_{{\bf
PGL}_r}^d,{\cal
D}^{k})$ for $g\ge 2$, {\it is not necessarily an integer} for $g=1$:  for
$d=0$ one
finds $\displaystyle N(1)=1+{(k+1)^{r-1}-1\over r^2}$, which is not an integer
unless
$r^2\mid k$.

\vfill\eject \centerline{ REFERENCES} \vglue15pt\baselineskip12.8pt
\def\num#1{\item{\hbox to\parindent{\enskip [#1]\hfill}}}
\parindent=1.3cm
\num{A-S} M.F.\ {\pc ATIYAH}, I.M.\ {\pc SINGER}: {\sl The index of elliptic
operators III}. Ann.\ of Math.\ {\bf 87}, 546-604 (1968).  \smallskip

\num{B} A.\ {\pc BEAUVILLE}: {\sl Conformal blocks, Fusion rings and the
Verlinde
formula.} Proc.\ of the Hirzebruch 65 Conf.\ on Algebraic Geometry, Israel
Math.\
Conf.\ Proc.\ {\bf 9}, 75-96 (1996).
\smallskip

\num{B-L} A.\ {\pc BEAUVILLE}, Y.\ {\pc LASZLO}: {\sl
Conformal blocks and generalized theta functions.} Comm.\
Math.\ Phys.\  {\bf 164}, 385-419 (1994).

\smallskip
\num{B-L-S} A.\ {\pc BEAUVILLE}, Y.\ {\pc LASZLO}, Ch.\ {\pc SORGER}: {\sl
The Picard group of the moduli of $G$\tx bundles on a
curve}. Preprint alg-geom/9608002.
\smallskip
\num{B-S} A. {\pc BERTRAM}, A. {\pc SZENES}: {\sl Hilbert polynomials of
moduli spaces of rank $2$ vector bundles II.} Topology {\bf 32}, 599-609
(1993). \smallskip
\num{Bo} N.\ {\pc BOURBAKI}: {\sl Alg\`ebre}, Chap.\ X (Alg\`ebre homologique).
Masson, Paris (1980).   \smallskip

 \num{D-N} J.M.\ {\pc DREZET}, M.S.\ {\pc NARASIMHAN}: {\sl Groupe de Picard
des vari\'et\'es de modules de fibr\'es semi-stables sur les courbes
alg\'ebriques.} Invent.\ math.\ {\bf 97}, 53-94 (1989).
 \smallskip

\num {F} G.\ {\pc FALTINGS}: {\sl A proof for the Verlinde formula.} J.\
Algebraic Geometry {\bf 3}, 347-374 (1994). \smallskip

\num{F-S-S}  J.\ {\pc FUCHS}, B.\ {\pc SCHELLEKENS}, Ch.\ {\pc SCHWEIGERT}:
{\sl From Dynkin diagram symmetries to fixed point structures}. Preprint
hep-th/9506135.
\smallskip
 \num{F-H} W.\ {\pc FULTON}, J.\ {\pc HARRIS}: {\sl Representation
theory}. GTM {\bf 129}, Springer-Verlag, New York Berlin Heidelberg (1991).
\smallskip

\num{L-S} Y.\ {\pc LASZLO}, Ch.\ {\pc SORGER}: {\sl The line bundles
on the moduli of parabolic $G$\tx bundles over curves and their
sections}. Annales de l'ENS, to appear; preprint alg-geom/9507002.
\smallskip

\num{M-S} G. {\pc MOORE}, N. {\pc SEIBERG}: {\sl Taming the conformal zoo}.
Phys.
Letters B {\bf 220}, 422-430 (1989). \smallskip

\num{N-R} M.S.\ {\pc NARASIMHAN}, S.\ {\pc RAMANAN}:
{\sl Generalized Prym varieties as fixed points}. J.\ of the Indian
Math.\ Soc.\ {\bf 39}, 1-19 (1975).
\smallskip
\num{P} T.\ {\pc PANTEV}: {\sl Comparison of generalized theta functions}. Duke
Math.\ J.\ {\bf 76}, 509-539 (1994).
\smallskip
\num{R} S.\ {\pc RAMANAN}:
{\sl The moduli spaces of vector bundles over an algebraic curve}. Math. Ann.
{\bf
200}, 69-84 (1973). \vskip1cm

\def\pc#1{\eightrm#1\sixrm}
\hfill\vtop{\eightrm\hbox to 5cm{\hfill Arnaud {\pc BEAUVILLE}\hfill}
 \hbox to 5cm{\hfill DMI -- \'Ecole Normale
Sup\'erieure\hfill} \hbox to 5cm{\hfill (URA 759 du CNRS)\hfill}
\hbox to 5cm{\hfill  45 rue d'Ulm\hfill}
\hbox to 5cm{\hfill F-75230 {\pc PARIS} Cedex 05\hfill}}

 \end